\documentclass[11pt,a4paper,reqno]{amsart}

\usepackage{mathptmx}
\usepackage[scaled=.90]{helvet}
\usepackage{courier}

\usepackage{graphicx,subfigure,parskip,amssymb,amsmath}
\newtheorem{theorem}{Theorem}
\newtheorem{lemma}[theorem]{Lemma}

\newcommand{\interior}{\operatorname{int}}
\newcommand{\D}{\mathbb D}

\begin{document}
\title{Approximating the modulus of an inner function}
\date{September 13, 2006. To appear in Pacific J.\ Math.}
\author{Geir Arne Hjelle}
\address{G.\ A.\ Hjelle\\Department of Mathematical Sciences\\Norwegian University of Science and Technology\\7491 Trondheim\\Norway}
\curraddr{G.\ A.\ Hjelle\\Department of Mathematics\\Washington University\\St.\ Louis, MO 63130\\USA}
\urladdr{http://www.math.wustl.edu/\textasciitilde hjelle/}
\email{hjelle@math.wustl.edu}
\author{Artur Nicolau}
\address{A.\ Nicolau\\Departament de Matem\`atiques\\Universitat Aut\`onoma de Barcelona\\08193 Bellaterra, Barcelona\\Spain}
\urladdr{http://mat.uab.es/\textasciitilde artur/}
\email{artur@mat.uab.es}
\subjclass[2000]{Primary 30D50, 30E10}
\keywords{Blaschke product, Carleson contour, discretization, inner function, interpolating Blaschke product, modulus}
\thanks{The first author is partially supported by grants from the Research Council of Norway, projects \#155060 and \#166395. The second author is partially supported by MTM2005-00544 and 2005SGR00774.}

\begin{abstract}
We show that the modulus of an inner function can be uniformly
approximated in the unit disk by the modulus of an interpolating
Blaschke product.
\end{abstract}
\maketitle

\section{Introduction}
Let $H^\infty$ be the algebra of bounded analytic functions in the
unit disk $\D$.  A function in $H^\infty$ is called inner if it has
radial limit of modulus one at almost every point of the unit
circle. A Blaschke product is an inner function of the form
\begin{equation*}
  B(z)=z^m \prod_{n=1}^\infty \frac{\bar z_n}{|z_n|} \;
    \frac{z_n-z}{1-\bar z_n z},
\end{equation*}
where $m$ is a non-negative integer and $\{z_n\}$ is a sequence of
points in $\D \setminus \{0\}$ satisfying the Blaschke condition
$\sum_n (1 - |z_n|) < \infty$. A classical result of O.\ Frostman
tells that for any inner function $f$, there exists an exceptional set
$E = E(f) \subset \D$ of logarithmic capacity zero such that the
M\"obius shift
\begin{equation*}
  \frac{f - \alpha}{1 - \bar \alpha f}
\end{equation*}
is a Blaschke product for any $\alpha \in \D \setminus E$.
See~\cite{Frostman35} or~\cite[p.\ 79]{Garnett81}. Hence any inner
function can be uniformly approximated by a Blaschke product.

A Blaschke product $B$ is called an interpolating Blaschke product if
its zero set $\{z_n\}$ form an interpolating sequence, that is, for
any bounded sequence of complex numbers $\{w_n\}$, there exists a
function $f \in H^\infty$ such that $f(z_n) = w_n$, $n = 1, 2, \ldots$. A
celebrated result by L.\ Carleson tells that this holds precisely when
the following two conditions are satisfied:
\begin{enumerate}
\item $\displaystyle{\inf_{n \neq m} \bigl| \frac{z_n - z_m}{1 - \bar
z_m z_n} \bigr| > 0}$,
\item there exists a constant $C$ such that $\sum_{z_n \in Q} (1 -
|z_n|) < C \ell(Q)$ for any Carleson square $Q$ of the form
\begin{equation}
\label{eq:Carlesonsquare}
  Q = \bigl\{ r e^{i \theta} \colon 0 < 1 - r < \ell(Q), \
      |\theta - \theta_0| < \pi \ell(Q) \bigr\}
\end{equation}
where $\theta_0 \in [0, 2\pi)$ and $0 < \ell(Q) < 1$.
\end{enumerate}
See~\cite{Carleson58} or~\cite[p.\ 287]{Garnett81}. Although the
interpolating Blaschke products comprise a small subset of all
Blaschke products, they play a central role in the theory of the
algebra $H^\infty$. See the last three chapters of~\cite{Garnett81}.

In~\cite{Marshall76} D.\ Marshall proved that any function $f\in
H^\infty$ can be uniformly approximated by finite linear combinations
of Blaschke products. That is, for any $\varepsilon >0$ there are
constants $c_1,\ldots, c_N$ and Blaschke products $B_1,\ldots, B_N$
such that
\begin{equation*}
  \bigl\| f - \sum_{i=1}^N c_i B_i \bigr\|_\infty < \varepsilon.
\end{equation*}
Here the $\infty$-norm is given by $\| g \|_\infty = \sup \{ | g(z)| :
z \in \D \}$.  This result was improved in~\cite{Garnett96} by showing
that one can take each of $B_1,\ldots, B_N$ to be an interpolating
Blaschke product. However the following problem remains open.

\begin{enumerate}
\item For any inner function $B$ and $\varepsilon > 0$, is there an
interpolating Blaschke product $I$ such that $\|B-I\|_\infty <
\varepsilon$?
\end{enumerate}

This question was posed in~\cite[p.\ 430]{Garnett81}, \cite[pp.\
268--269]{Havin94}, \cite{Jones81} and~\cite[p.\ 202]{Nikolskii86}.
The purpose of this note is to provide a positive answer if one
restricts attention to the modulus.

\begin{theorem}
\label{thm:maintheorem}
Let $B$ be an inner function and $\varepsilon >0$. Then there
exists an interpolating Blaschke product $I$ such that
\begin{equation*}
  \bigl| |B(z)|- |I(z)| \bigr| < \varepsilon
\end{equation*}
for all $z \in \D$.
\end{theorem}

The proof may be described as follows. The first step consists of
constructing a system $\Gamma = \bigcup_i \Gamma_i$ of disjoint closed
curves $\Gamma_i \subset \D$ such that arclength of $\Gamma$ is a
Carleson measure, and verifying that
\begin{enumerate}
\item[(a)] $|B(z)|$ is uniformly small on hyperbolic disks of fixed
radius centered at points of~$\Gamma$,
\item[(b)] in any hyperbolic disk of fixed radius centered at a point
outside the union of the interiors of $\Gamma_i$, $\bigcup_i \interior
\Gamma_i$, there is a point $z$ where $|B(z)|$ is not small.
\end{enumerate}
Write $B = B_1 \cdot B_2$ where $B_1$ is the Blaschke product formed
with the zeros of $B$ which are in $\bigcup_i \interior \Gamma_i$.
Statement (b) gives that $B_2$ is a finite product of interpolating
Blaschke products. Since D.\ Marshall and A.\ Stray proved
in~\cite{Marshall96} that any finite product of interpolating Blaschke
products may be approximated by a single interpolating Blaschke
product, the relevant zeros of $B$ lie in $\bigcup_i \interior
\Gamma_i$, that is, are those of $B_1$. The construction of $\Gamma$
is a variation of the original Corona construction introduced by L.\
Carleson. See~\cite{Carleson62} or~\cite[pp.\ 342--347]{Garnett81}.

Next, for each $i = 1, 2, \ldots$, let $\mu_i$ be the sum of harmonic
measures in $\interior \Gamma_i$ from the zeros of $B_1$ contained in
$\interior \Gamma_i$. Then the mass $\mu_i(\Gamma_i)$ is the total
number of zeros of $B_1$ contained in $\interior \Gamma_i$. The second
step consists of splitting $\Gamma_i = \bigcup_k \Gamma_{i,k}$, into
pieces $\Gamma_{i,k}$ with $\mu_i(\Gamma_{i,k}) = 1$, $k = 1,
2,\ldots$ and choosing points $\xi_{i,k}\in \Gamma_{i,k}$ which match
a certain moment of the measure $\mu_i$ on $\Gamma_{i,k}$.  This
choice may be compared with~\cite{Lyubarskii01} where a related
discretization argument is performed in a different context. Let $I_1$
be the Blaschke product with zeros $\xi_{i,k}$, $i, k = 1, 2, \ldots$.
Finally the last step of the proof is to use $(b)$ above to show that
$I_1$ is an interpolating Blaschke product and to use the location of
$\{ \xi_{i,k} \}$, as well as (a) above, to show that $|I_1(z) \cdot
B_2(z)| $ approximates $|B(z)|$.

Besides the individual problem mentioned above, some questions
concerning approximation by arguments of interpolating Blaschke
products remain open. Let $B$ be an inner function.
\begin{enumerate}
\addtocounter{enumi}{1} 
\item Given $\varepsilon >0$, is there an interpolating Blaschke
product $I$ such that
\begin{equation*}
  \|\operatorname{Arg} B -\operatorname{Arg} I\|_{\operatorname{BMO}
    (\partial \D)} < \varepsilon ?
\end{equation*}
\item Is there an interpolating Blaschke product $I$ such that
$\operatorname{Arg} B - \operatorname{Arg} I = \widetilde v$ where $v
\in L^\infty (\partial \D)$?
\item Is there an interpolating Blaschke product $I$ such that
$\operatorname{Arg} B -\operatorname{Arg} I = u + \tilde v$ where $u,
v \in L^\infty (\partial \D)$ and $\| u\|_\infty < \frac\pi2$?
\end{enumerate}

It is clear that a positive answer to any of these problems would lead
to a positive answer to the next one. Moreover a positive answer to
Problem 2 would imply the main result of this note.  Problem 4 was
posed by N.\ K.\ Nikol$'$ski{\u\i} in~\cite{Havin94}
and~\cite{Nikolskii86} in connection to Toeplitz operators and
complete interpolating sequences in model spaces. Problem 3 and
Problem 4 have been discussed in the nice monograph by K.\
Seip~\cite[p.\ 92]{Seip04}.

We are indebted to Arne Stray for his valuable comments on an earlier
version of this paper. Part of this work was done while the first
author was visiting Universitat Aut\`onoma de Barcelona, and while the
second author was visiting IMUB at Universitat de Barcelona. It is a
pleasure to thank both institutions for their support.

\section{Construction of the contour}
The hyperbolic distance between two points $z, w \in \D$ is
\begin{equation*}
  \beta(z, w) = \tfrac12 \log \tfrac{1 + \rho(z, w)}{1 - \rho(z, w)}
\end{equation*}
where $\rho(z, w)$ is the pseudohyperbolic distance,
\begin{equation*}
  \rho(z, w) = \bigl| \frac{z - w}{1 - \bar w z} \bigr| .
\end{equation*}

Recall that a positive measure $\mu$ in the unit disk is called a
Carleson measure if there exists a constant $M=M(\mu)>0$ such that
$\mu (Q)\leq M\ell (Q)$ for any Carleson square of the
form~\eqref{eq:Carlesonsquare}. The infimum of the constants $M$
verifying the inequality above is called the Carleson norm of the
measure $\mu$ and it is denoted by $\|\mu\|_C$.

The main result of this section is a variant of the classical
construction of the Carleson contour introduced by L. Carleson in his
original proof of the Corona Theorem.  See~\cite{Carleson62}
or~\cite[pp.\ 342--347]{Garnett81}.

\begin{lemma}
\label{lem:contour}
Let $B \in H^\infty$ with $\|B\|_\infty = 1$. Let $0 < \varepsilon <1$
and $K > 0$ be fixed constants. Then, there exists a constant $\delta
= \delta(\varepsilon, K) > 0$ and a system $\Gamma = \bigcup \Gamma_i$
of disjoint closed curves $\Gamma_i$ contained in $\D$ such that
\begin{enumerate}
\item[(a)] if $\inf_i \beta (z, \interior \Gamma_i)\leq K$, one has
$|B(z)|\leq \varepsilon$,
\item[(b)] if $z\notin \bigcup \interior \Gamma_i$, one has
\begin{equation*}
  \sup \{|B(w)| \colon \beta(w, z) \leq K + 14 \} > \delta ,
\end{equation*}
\item[(c)] arclength on $\Gamma$, $\mathrm ds_{|\Gamma}$, is a Carleson
measure and $\|\mathrm ds_{|\Gamma}\|_C \leq 68$.
\end{enumerate}
\end{lemma}

\begin{proof}
The proof is essentially contained in the recent
paper~\cite{Nicolau05}, but we sketch it for the convenience of the
reader. Given a set $E \subset \D$, let $\Omega_K(E)$ denote the set
of points that are at most at hyperbolic distance $K$ from the set
$E$, that is,
\begin{equation*}
  \Omega_K(E) = \bigl\{ z \colon \inf_{w \in E} \beta(z, w) \leq K \bigr\} .
\end{equation*}
Consider dyadic Carleson squares of the form
\begin{equation*}
  Q_{n,j} = \bigl\{ r e^{i \theta} \colon 1 - 2^{-n} < r < 1, \ 2\pi j 2^{-n}
    < \theta < 2\pi (j+1) 2^{-n} \bigr\},
\end{equation*}
for $j = 0, 1, \ldots, 2^n-1$ and $n = 1, 2, \ldots$, and their top
halves $T(Q_{n,j})= \{r e^{i \theta}\in Q_{n,j} \colon r < 1-2^{-n-1}
\}$. Let $0 < \delta < \varepsilon$ be a constant to be fixed later.
A dyadic Carleson square $Q$ will be called good if
\begin{equation*}
  \sup \bigl\{ |B(z)| \colon z \in \Omega_K \bigl( T(Q) \bigr) \bigr\}
    > \varepsilon .
\end{equation*}
The collection of good dyadic Carleson squares will be denoted by
$\{Q_j^G \colon j = 1, 2, \ldots\}$. A dyadic Carleson square $Q$ will be
called bad if
\begin{equation*}
  \sup \bigl\{|B(z)| \colon z \in \Omega_K \bigl( T(Q) \bigr) \bigr\}
    < \delta .
\end{equation*}
We denote the collection of bad dyadic Carleson squares by $\{Q_j^B
\colon j = 1, 2, \ldots\}$. The construction goes as follows.
\begin{figure}[bt]
\begin{center}
  {\includegraphics{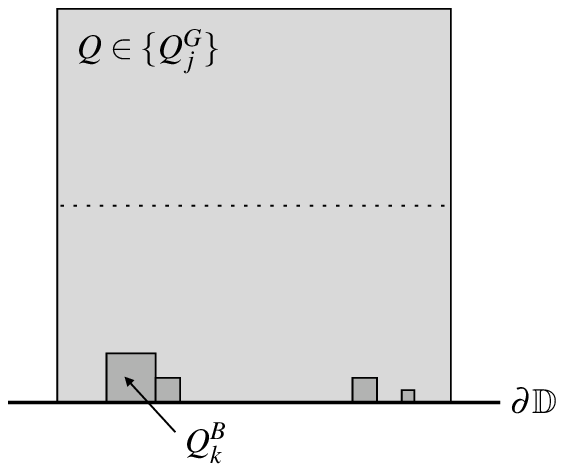}} \hfil
  {\includegraphics{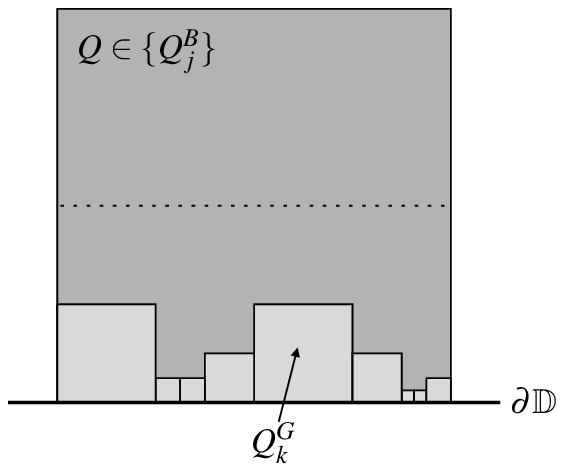}}
\caption{Choosing good and bad squares for constructing the contour}
\label{fig:contour_gb}
\end{center}
\end{figure}

1. For each good dyadic Carleson square $Q = Q_j^G$, we choose the
maximal bad dyadic Carleson squares $Q_k^B$ contained in $Q$. The main
estimate in the construction is
\begin{equation}
\label{eq:scaling}
  \sum_{Q_k^B \subset Q} \ell(Q_k^B) \leq \frac12 \ell(Q).
\end{equation}
Since $|B(z)| < \delta$ if $z \in T(Q_k^B)$, while $|B(z)| >
\varepsilon$ for some $z \in \Omega_K(T(Q))$, taking $\delta =
\delta(\varepsilon, K)$ sufficiently small, standard arguments lead
to~\eqref{eq:scaling}. See Lemma~2.1 of~\cite{Nicolau05} for
details.

2. For each bad dyadic Carleson square $Q = Q_j^B$, we choose the
maximal good dyadic Carleson squares $Q_k^G$ contained in $Q$.  This
family is denoted by $G(Q) = \{ Q_k^G \colon k = 1, 2, \ldots\}$.
\begin{figure}[bt]
\begin{center}
\includegraphics{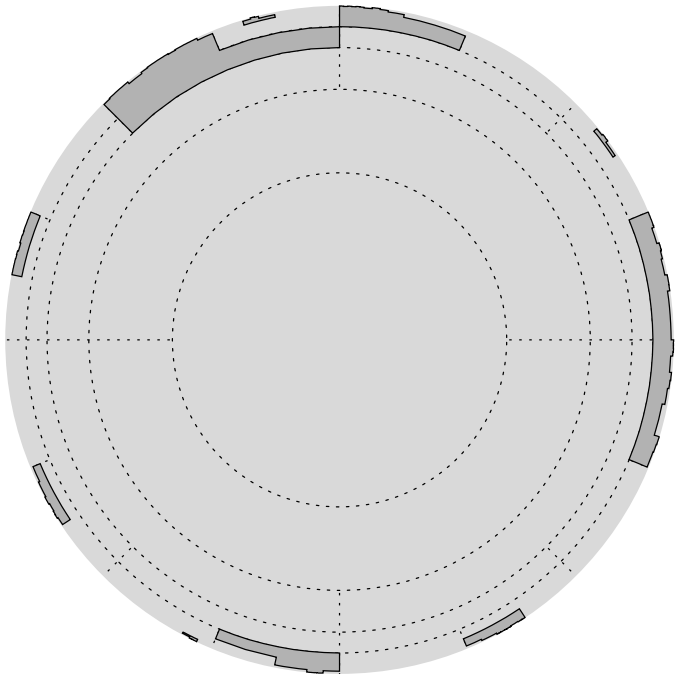}
\caption{The unit disk, some dyadic Carleson contours and an example of a contour.}
\label{fig:contour_circle}
\end{center}
\end{figure}

So, from each good dyadic Carleson square we move to bad ones
fulfilling the estimate~\eqref{eq:scaling} and from each bad one we
again move to good ones. See Figure~\ref{fig:contour_gb}. Now for each
bad square $Q = Q_j^B$, let $R(Q)$ be the region
\begin{equation*}
  R(Q) = Q \setminus \overline{\bigcup_{G(Q)} Q_k^G}
\end{equation*}
and let $R$ be the open set
\begin{equation*}
  R = \bigcup_j R(Q_j^B) .
\end{equation*}
Finally, decompose $R$ into its connected components $R_i$ and denote
$\Gamma_i = \partial R_i$, $i = 1, 2, \ldots$. Observe that each
$\Gamma_i$ consists of pieces of boundaries of dyadic Carleson
squares. See Figure~\ref{fig:contour_circle}. By construction if $z
\in R$ we have
\begin{equation*}
  \sup \bigl\{|B(w)| \colon \beta(w, z) \leq K \bigr\} \leq \varepsilon
\end{equation*}
and hence part $(a)$ in the statement follows. Similarly, if $z \notin
R$, the point $z$ is not in the top part of a bad dyadic Carleson
square. As the hyperbolic diameter of a top part of a Carleson square
is uniformly bounded, say by 14, we deduce that there exists $w \in
\D$ with $\beta(z,w) \leq K + 14$ such that $|B(w)| > \delta$. Hence
part $(b)$ in the statement follows.  Since the length of $\partial
R(Q)$ is bounded by $17 \ell(Q)$, the scaling~\eqref{eq:scaling} shows
that for any bad dyadic square $Q$, one has
\begin{equation*}
  \sum_{Q_j^B \subsetneq Q} |\partial R(Q_j^B)|\leq 17 \ell (Q).
\end{equation*}
Then easy geometric considerations show that arclength on
$\bigcup\Gamma_i$ is a Carleson measure and its Carleson norm is
smaller than $68$.
\end{proof}

\section{Construction of the interpolating Blaschke product}
We now use Lemma~\ref{lem:contour} to construct a contour
$\Gamma$. Note that by Frostman's Theorem we can assume that $B$ is a
Blaschke product. Given $\varepsilon > 0$, let $N$ be a big constant
dependent on $\varepsilon$ to be fixed later. Apply
Lemma~\ref{lem:contour} with $\frac\varepsilon2$ and $2N$ instead of
$\varepsilon$ and $K$ to obtain $\Gamma$ and $\delta > 0$ such that
\begin{enumerate}
\item[(a)] $|B(z)| < \frac\varepsilon2$ if $\beta(z, \interior \Gamma)
\leq 2N$,
\item[(b)] $\sup \{ |B(w)| \colon \beta(w, z) \leq 2N + 14 \} > \delta$
if $z \not\in \interior \Gamma$,
\item[(c)] arclength on $\Gamma$ is a Carleson measure with Carleson
norm $\|\mathrm ds_{|\Gamma}\|_C \leq 68$.
\end{enumerate}
With the contour $\Gamma$ in place, we want to construct the
interpolating Blaschke product $I$. Split $B$ into two Blaschke
products $B_1$ and $B_2$. That is $B = B_1 \cdot B_2$, where $B_1$ is
formed with the zeros $\{z_n\}$ of $B$ which are inside $\interior
\Gamma$ and at hyperbolic distance more than 1 from the contour
$\Gamma$. Now for each zero $z$ of $B_2$, part (b) provides a point $w
\in \D$, $\beta(w, z) \leq 2N + 15$ such that $|B_2(w)| \geq |B(w)| >
\delta$. This implies that $B_2$ is a finite product of interpolating
Blaschke products (see Theorem 2.2 of~\cite{Mortini04}).

\begin{figure}[bt]
\begin{center}
\includegraphics{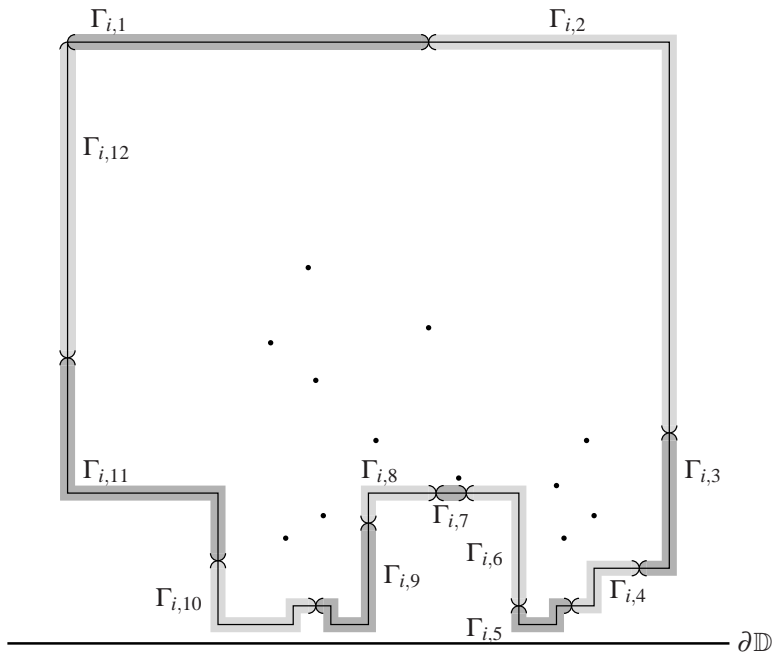}
\caption{Each component $\Gamma_i$ of the contour is split into arcs $\Gamma_{i,k}$ such that the $\mu$-measure of each arc is $1$.}
\label{fig:contoursplit}
\end{center}
\end{figure}

Hence the dangerous part of $B$ will be $B_1$ which has all its zeros
contained deeply inside the contour $\Gamma$. We want to mimic the
behavior of $|B_1|$ by constructing a Blaschke product $I_1$ with
zeros on $\Gamma$. To this end, for each component $\Gamma_i$ of the
contour we consider the measure
\begin{equation*}
  \mathrm d\mu_i(\xi) = \sum_{\substack{z_n \in \interior \Gamma_i\\
    \beta(z_n, \Gamma_i) > 1}}
    \omega(z_n, \xi; \interior \Gamma_i)
\end{equation*}
defined for $\xi \in \Gamma_i$. Here $\omega(z, \xi; \Omega)$ denotes
the harmonic measure from the point $z \in \Omega$ in the domain
$\Omega \subseteq \D$. Clearly $\mu_i(\Gamma_i)$ will be equal to the
number of zeros $z_n$ of $B_1$ inside $\Gamma_i$. Next we split
$\Gamma_i$ into disjoint arcs $\Gamma_{i,k}$ such that
$\mu_i(\Gamma_{i,k}) = 1$ for each $k$. This is illustrated in
Figure~\ref{fig:contoursplit}. On each such arc we locate one zero
$\xi_{i,k}$ of $I_1$ such that
\begin{equation}
\label{eq:placing_xi-i,k}
  1 - |\xi_{i,k}|^2 = \int_{\Gamma_{i,k}}
      \bigl( 1 - |\xi|^2 \bigr) \, \mathrm d\mu_i(\xi) .
\end{equation}
This will in general not determine the points $\xi_{i,k}$
uniquely. However, there seems to be a lot of freedom for placing the
zeros of $I_1$ in this construction, and the
condition~\eqref{eq:placing_xi-i,k} will be sufficient for our
purposes.

Let $I_1$ be the Blaschke product with the zeros $\xi_{i,k}$, and
factor $I_1 = I_1^o \cdot I_1^e$ where $I_1^o$ is the Blaschke product
with zeros $\xi_{i,k}$ with $k$ odd, while $I_1^e$ is the Blaschke
product with zeros $\xi_{i,k}$ with $k$ even. In
Figure~\ref{fig:contoursplit}, $I_1^o$ has its zeros placed in the
dark arcs, while the zeros of $I_1^e$ are placed in the light arcs. We
claim that both $I_1^o$ and $I_1^e$ are interpolating Blaschke
products, and hence $I_1$ can be approximated by an interpolating
Blaschke product~\cite{Marshall96}. To show this claim we will observe
that their zero sets satisfy the two conditions of Carleson's
theorem~\cite{Carleson58}, stated in the introduction.

In this case, property (2) follows from the fact that arclength is a
Carleson measure on $\Gamma$, while the first property follows from
the following lemma and the geometry of the contour.
\begin{lemma}
The hyperbolic length, $\ell_\beta(\Gamma_{i,k})$, of $\Gamma_{i,k}$
is bounded from below,
\begin{equation*}
  \ell_\beta(\Gamma_{i,k}) \geq \delta^{2 e^{2(2N + 14)}} .
\end{equation*}
\end{lemma}

\begin{proof}
We first show that for any point $w \in \Gamma$, $|B_1(w)|$ is bounded
from below by some constant depending only on $\delta$ and $N$. To see
this, recall that there is a point $\zeta$ such that $\beta(\zeta, w)
\leq 2N + 14$ and $|B_1(\zeta)| \geq |B(\zeta)| > \delta$. Consider
\begin{equation*}
  \log \bigl| B_1(w) \bigr|^{-1}
    = \sum_{z_n} \log \rho(w, z_n)^{-1} ,
\end{equation*}
where the sum is taken over all zeros $z_n$ of $B_1$. As $w$ is
separated from the zeros of~$B_1$,
\begin{equation*}
  \log \rho(w, z_n)^{-1} \leq 1 - \rho(w, z_n)^2 .
\end{equation*}
Furthermore,
\begin{equation*}
  \rho(w, z_n)
    \geq \frac{\rho(z_n, \zeta) - \rho(\zeta, w)}
      {1 - \rho(z_n, \zeta) \rho(\zeta, w)}
    \geq \frac{\rho(z_n, \zeta) - C}{1 - C \rho(z_n, \zeta)},
\end{equation*}
where $C = \frac{e^{2(2N+14)} - 1}{e^{2(2N+14)} + 1} < 1$. Hence
\begin{align*}
  \log \rho(w, z_n)^{-1}
    &\leq \frac{\bigl( 1 - \rho(z_n, \zeta)^2 \bigr) \bigl(1 - C^2 \bigr)}
      {\bigl( 1 - C \rho(z_n, \zeta) \bigr)^2}
    \leq \frac{1 + C}{1 - C} \bigl( 1 - \rho(z_n, \zeta)^2 \bigr) \\
    &\leq 2 \frac{1 + C}{1 - C} \log \rho(z_n, \zeta)^{-1}
    = 2 e^{2(2N + 14)} \log \rho(z_n, \zeta)^{-1} ,
\end{align*}
and we see that $|B_1(w)| \geq \delta^{2 e^{2(2N + 14)}}$.

Intuitively, this lower bound for the values of $|B_1|$ should imply
that the arcs $\Gamma_{i,k}$ can not be too short hyperbolically. To
make this observation rigorous we argue as follows. Using that the
harmonic measure $\omega$ is positive and harmonic, we have that for
any $z \in \interior \Gamma_i$,
\begin{equation*}
  \omega(z, \Gamma_{i,k}; \interior \Gamma_i)
    \leq \omega(z, \Gamma_{i,k}; \D \setminus \Gamma_{i,k})
    \leq \frac{\displaystyle{\int_{\Gamma_{i,k}} \log \bigl|
      \frac{z - w}{1 - \bar w z} \bigr|^{-1} \, \frac{|\mathrm dw|}
      {1 - |w|^2}}}{\displaystyle{\min_{z \in \Gamma_{i,k}}
      \int_{\Gamma_{i,k}} \log \bigl| \frac{z - w}{1 - \bar w z}
      \bigr|^{-1} \, \frac{|\mathrm dw|}{1 - |w|^2}}}
\end{equation*}
and
\begin{align*}
  1 = \mu_i(\Gamma_{i,k})
    &= \sum_{z_n \in \interior \Gamma_i} \omega(z_n, \Gamma_{i,k};
      \interior \Gamma_i) \\
    &\leq \frac1{C_{i,k}} \int_{\Gamma_{i,k}} \log
      \biggl( \prod_{z_n \in \interior \Gamma_i} \bigl|
      \frac{z_n - w}{1 - \bar w z_n} \bigr|^{-1} \biggr) \,
      \frac{|\mathrm dw|}{1 - |w|^2}
\end{align*}
where $C_{i,k} = \min_{z \in \Gamma_{i,k}} \int_{\Gamma_{i,k}} \log
\bigl| \frac{z - w}{1 - \bar w z} \bigr|^{-1} \, \frac{|\mathrm dw|}{1
- |w|^2}$ is a constant dependent on $\Gamma_{i,k}$. Let $B_{1,i}$
denote the Blaschke product with the zeros of $B_1$ that fall inside
the component $\Gamma_i$. Then for $w \in \Gamma_i$,
\begin{equation*}
  \log \biggl( \prod_{z_n \in \interior \Gamma_i} \bigl|
      \frac{z_n - w}{1 - \bar w z_n} \bigr|^{-1} \biggr)
    = \log |B_{1,i}(w)|^{-1}
    \leq \log |B_1(w)|^{-1}
    \leq 2 e^{2(2N + 14)} \log \delta^{-1} .
\end{equation*}
Thus
\begin{equation*}
  1 \leq \frac1{C_{i,k}} 2 e^{2(2N + 14)} \log \delta^{-1}
      \int_{\Gamma_{i,k}} \, \frac{|\mathrm dw|}{1 - |w|^2}
    = \frac1{C_{i,k}} 2 e^{2(2N + 14)} \log \delta^{-1}
      \ell_\beta(\Gamma_{i,k})
\end{equation*}
such that
\begin{equation*}
  \ell_\beta(\Gamma_{i,k}) \geq \frac{C_{i,k}}{2 e^{2(2N + 14)}
    \log \delta^{-1}} .
\end{equation*}

To estimate $C_{i,k}$ we use the substitution $\xi = \varphi_z(w) =
\frac{z - w}{1 - \bar w z}$ and the conformal invariance of the
hyperbolic metric. A calculation then gives that
\begin{equation*}
  C_{i,k} \geq \log \bigl( \tanh(\ell_\beta(\Gamma_{i,k})) \bigr)
    \ell_\beta(\Gamma_{i,k}) ,
\end{equation*}
which implies the desired bound, $\ell_\beta(\Gamma_{i,k}) \geq
\delta^{2 e^{2(2N + 14)}}$.
\end{proof}

\section{Proof of the approximation}
In this section we will show that the constructed function, $I = I_1
\cdot B_2$, approximates the given Blaschke product uniformly in
modulus.  We first claim that it suffices to prove
Theorem~\ref{thm:maintheorem} for points $z \in \D$ far away from the
contour. Indeed, assume that we can prove that
\begin{equation}
\label{eq:faraway_approx}
  \bigl| |B_1(z)| - |I_1(z)| \bigr| < \tfrac\varepsilon2
\end{equation}
for all $z$ such that $\beta(z, \interior \Gamma) \geq 2N$, where $N$
is as in the construction of the contour. Then for points $z$ with
$\beta(z, \interior \Gamma) = 2N$
\begin{align*}
  |I(z)| &= \bigl(|I_1(z)| - |B_1(z)| + |B_1(z)|\bigr) |B_2(z)| \\
    &\leq \bigl| |B_1(z)| - |I_1(z)| \bigr| + |B(z)|
    < \tfrac\varepsilon2 + \tfrac\varepsilon2 = \varepsilon .
\end{align*}
By the maximum principle $|I(z)| < \varepsilon$ for all $z \in
\Omega_{2N}(\interior \Gamma)$ as well. Hence
\begin{equation*}
  \bigl| |B(z)| - |I(z)| \bigr|
    = \bigl| |B_1(z)| - |I_1(z)| \bigr| |B_2(z)|
    < \begin{cases}
        \frac\varepsilon2 & \text{if } \beta(z, \interior \Gamma) \geq 2N , \\
        \varepsilon & \text{if } \beta(z, \interior \Gamma) < 2N .
      \end{cases}
\end{equation*}
So Theorem~\ref{thm:maintheorem} follows
from~\eqref{eq:faraway_approx}.

The rest of the paper will be dedicated to prove
that~\eqref{eq:faraway_approx} holds. Fix a point $z$ such that
$\beta(z, \interior \Gamma) \geq 2N$. We will consider the logarithm
of $|B_1|$. As all the zeros of $B_1$ lie inside the contour $\Gamma$,
$\log \bigl| \frac{z - z_n}{1 - \bar z_n z} \bigr|$ is harmonic inside
$\Gamma$ as a function of $z_n$. Hence
\begin{equation*}
  \log |B_1(z)| = \sum_j \log \bigl| \frac{z - z_n}{1 - \bar z_n z} \bigr|
    = \int_\Gamma \log \bigl| \frac{z - \xi}{1 - \bar \xi z} \bigr| \,
      \mathrm d\mu(\xi) ,
\end{equation*}
where $\mathrm d\mu = \sum_i \mathrm d\mu_i$. As the $\mu$-measure of
each arc $\Gamma_{i,k}$ is $1$, we have
\begin{align}
  \log |B_1(z)| - \log |I_1(z)|
    &= \int_\Gamma \log \bigl| \frac{z - \xi}{1 - \bar \xi z} \bigr| \,
      \mathrm d\mu(\xi) - \sum_{i,k} \log \bigl| \frac{z - \xi_{i,k}}
      {1 - \bar \xi_{i,k} z} \bigr| \notag \\
    &= \sum_{i,k} \int_{\Gamma_{i,k}} \biggl( \log \bigl| \frac{z - \xi}
      {1 - \bar \xi z} \bigr| - \log \bigl| \frac{z - \xi_{i,k}}
      {1 - \bar \xi_{i,k} z} \bigr| \biggr) \, \mathrm d\mu(\xi) \notag \\
    \label{eq:decomposable_sum}
    &= \sum_{i,k} \int_{\Gamma_{i,k}} \log \frac{\rho(z, \xi)}
      {\rho(z, \xi_{i,k})} \, \mathrm d\mu(\xi)
    \overset{\text{def}}{=} \sum_{i,k} H_{i,k}(z) .
\end{align}
To estimate this sum we consider different types of arcs. By $Q_z$ we
denote the Carleson square with $z$ as the midpoint on the
top-side. We say that an arc $\Gamma_{i,k}$ is in the class $\mathcal
B$ if $\Gamma_{i,k} \subset 2^N Q_z$. Note that since $\beta(z,
\interior \Gamma) \geq 2N$, this implies that such an arc lies very
close to the boundary. The rest of the arcs we split into short and
long arcs. For $n \geq N + 1$ define
\begin{equation*}
  \mathcal S_n = \bigl\{ \Gamma_{i,k} \colon \ell_\beta(\Gamma_{i,k}) < 1, \
    \Gamma_{i,k} \subset 2^n Q_z \bigr\} \setminus \bigl(\mathcal B \cup
    \bigcup_{i < n} \mathcal S_i \bigr)
\end{equation*}
and
\begin{equation*}
  \mathcal L_n = \bigl\{ \Gamma_{i,k} \colon \ell_\beta(\Gamma_{i,k}) \geq 1, \
    \Gamma_{i,k} \subset 2^n Q_z \bigr\} \setminus \bigl(\mathcal B \cup
    \bigcup_{i < n} \mathcal L_i \bigr) .
\end{equation*}
\begin{figure}[bt]
\begin{center}
\includegraphics{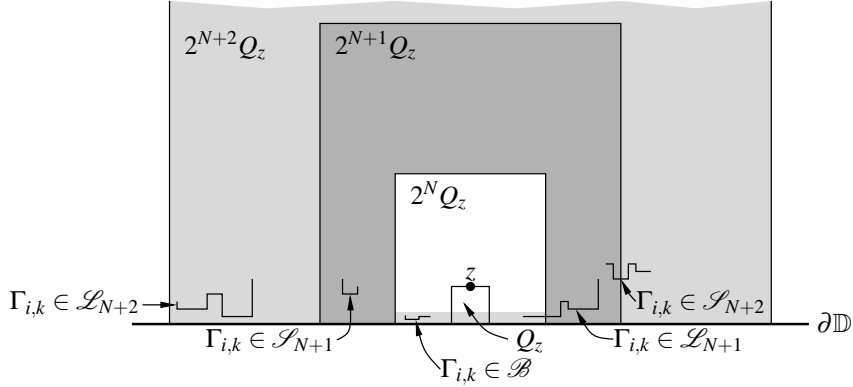}
\caption{We divide the arcs $\Gamma_{i,k}$ into classes denoted $\mathcal B$, $\mathcal S_n$ and $\mathcal L_n$}
\label{fig:qz_bst}
\end{center}
\end{figure}

Consult Figure~\ref{fig:qz_bst} for some examples of this
classification. This partition is such that each arc $\Gamma_{i,k}$
belongs to one and only one of the classes $\mathcal B$, $\mathcal
S_n$ and $\mathcal L_n$, $n \geq N + 1$. Hence we may decompose the
sum~\eqref{eq:decomposable_sum} as follows
\begin{equation*}
  \sum_{i,k} H_{i,k}(z) = \sum_{\Gamma_{i,k} \in \mathcal B} H_{i,k}(z)
    + \sum_{n = N+1}^\infty \biggl(
      \sum_{\Gamma_{i,k} \in \mathcal S_n} H_{i,k}(z) +
      \sum_{\Gamma_{i,k} \in \mathcal L_n} H_{i,k}(z) \biggr) .
\end{equation*}
Our goal is to show that the absolute value of the left hand side is
small. To accomplish this we will show that each of the terms
\begin{equation*}
  \biggl| \sum_{\Gamma_{i,k} \in \mathcal B} H_{i,k}(z) \biggr|, \quad
  \biggl| \sum_{n = N+1}^\infty \sum_{\Gamma_{i,k} \in \mathcal S_n}
    H_{i,k}(z) \biggr| \quad \text{and} \quad
  \biggl| \sum_{n = N+1}^\infty \sum_{\Gamma_{i,k} \in \mathcal L_n}
    H_{i,k}(z) \biggr|
\end{equation*}
are small.

Let us begin with the boundary arcs $\Gamma_{i,k} \in \mathcal
B$. Using that $\log (1 - t) = -t + \mathcal O(t^2)$ we get
\begin{multline*}
  \sum_{\Gamma_{i,k} \in \mathcal B} \int_{\Gamma_{i,k}} \log
      \frac{\rho(z, \xi)}{\rho(z, \xi_{i,k})} \, \mathrm d\mu(\xi) \\
    = -\frac12 \sum_{\Gamma_{i,k} \in \mathcal B} \int_{\Gamma_{i,k}}
      \biggl( 1 - \frac{\rho(z, \xi)^2}{\rho(z, \xi_{i,k})^2} +
      \mathcal O\biggl( \bigl( 1 - \frac{\rho(z, \xi)^2}
      {\rho(z, \xi_{i,k})^2} \bigr)^2 \biggr) \biggr) \, \mathrm d\mu(\xi) .
\end{multline*}
Taking absolute values,
\begin{multline}
\label{eq:boundaryarcs}
  \biggl| \sum_{\Gamma_{i,k} \in \mathcal B} \int_{\Gamma_{i,k}} \log
      \frac{\rho(z, \xi)}{\rho(z, \xi_{i,k})} \, \mathrm d\mu(\xi) \biggr|
    \leq \frac12 \biggl| \sum_{\Gamma_{i,k} \in \mathcal B}
      \int_{\Gamma_{i,k}} 1 - \frac{\rho(z, \xi)^2}{\rho(z, \xi_{i,k})^2}
      \, \mathrm d\mu(\xi) \biggr| \\ + \frac12 \biggl| \sum_{\Gamma_{i,k}
      \in \mathcal B} \int_{\Gamma_{i,k}} \mathcal O\biggl( \bigl( 1 -
      \frac{\rho(z, \xi)^2}{\rho(z, \xi_{i,k})^2} \bigr)^2 \biggr) \,
      \mathrm d\mu(\xi) \biggr|
    \overset{\text{def}}{=} E_{\mathcal B, 1} + E_{\mathcal B, 2} ,
\end{multline}
where we define $E_{\mathcal B, 1}$ and $E_{\mathcal B, 2}$ for
convenience.  At first we focus on the first term, $E_{\mathcal
B, 1}$, of this expression. Note that as $z$ is far away from
$\xi_{i,k} \in \Gamma_i$, $\rho(z, \xi_{i,k})^{-2}$ is bounded, say
$\rho(z, \xi_{i,k})^{-2} \leq 2$. By expanding $1 - \rho(z, \xi)^2$
and $1 - \rho(z, \xi_{i,k})^2$, we can write
\begin{align}
\label{eq:zerocancelling}
  E_{\mathcal B, 1} &\leq \sum_{\Gamma_{i,k} \in \mathcal B} \biggl|
      \int_{\Gamma_{i,k}} \bigl( 1 - |z|^2 \bigr) \biggl( \frac{1 - |\xi|^2}
      {|1 - \bar \xi z|^2} - \frac{1 - |\xi_{i,k}|^2}{|1 - \bar \xi_{i,k} z|^2}
      \biggr) \, \mathrm d\mu(\xi) \biggr| \\
    &= \sum_{\Gamma_{i,k} \in \mathcal B} \biggl| \int_{\Gamma_{i,k}}
      \bigl( 1 - |z|^2 \bigr) \biggl( \frac{1 - |\xi|^2}{|1 - \bar \xi z|^2}
      - \frac{1 - |\xi|^2}{|1 - \bar \xi_{i,k} z|^2} + \frac{|\xi_{i,k}|^2 - 
      |\xi|^2}{|1 - \bar \xi_{i,k} z|^2} \biggr) \, \mathrm d\mu(\xi) \biggr| .
      \notag
\end{align}
By the placement,~\eqref{eq:placing_xi-i,k}, of the zeros $\xi_{i,k}$,
the integral of the last term is zero. We now move the modulus under
the integral to get
\begin{equation}
\label{eq:movingmodulus}
  E_{\mathcal B, 1} \leq \bigl( 1 - |z|^2 \bigr) \sum_{\Gamma_{i,k} \in
      \mathcal B} \int_{\Gamma_{i,k}} \bigl( 1 - |\xi|^2 \bigr) \biggl|
      \frac1{|1 - \bar \xi z|^2} - \frac1{|1 - \bar \xi_{i,k} z|^2}
      \biggr| \, \mathrm d\mu(\xi) .
\end{equation}
Because $\xi$ and $\xi_{i,k}$ should be close to each other in some
sense, compared to $z$, we suspect some cancellation. Therefore we use
the estimate
\begin{equation}
\label{eq:cancellation}
  \biggl| \frac1{|1 - \bar \xi z|^2} - \frac1{|1 - \bar \xi_{i,k} z|^2} \biggr|
    \leq \frac{2 |\xi - \xi_{i,k}|}{\bigl( 1 - |z| \bigr)^3}
\end{equation}
and the more trivial inequalities $|\xi - \xi_{i,k}| \leq
\ell(\Gamma_{i,k})$ and $1 - |z|^2 \leq 2 (1 - |z|)$ to obtain
\begin{equation*}
  E_{\mathcal B, 1} \leq 2^3 \cdot \bigl( 1 - |z| \bigr)^{-2}
      \sum_{\Gamma_{i,k} \in \mathcal B} \ell(\Gamma_{i,k})
      \int_{\Gamma_{i,k}} \bigl( 1 - |\xi| \bigr) \, \mathrm d\mu(\xi) .
\end{equation*}
All the arcs $\Gamma_{i,k} \in \mathcal B$ are contained in a
rectangle at the boundary with height $2^{-2N} (1 - |z|)$ and width
$2^N (1 - |z|)$. Using that $1 - |\xi| \leq 2^{-2N} (1 - |z|)$ and
that the arclength $\mathrm ds_{|\Gamma}$ is a Carleson measure, we
then get
\begin{equation*}
  E_{\mathcal B, 1} \leq 2^3 \cdot \|\mathrm ds_{|\Gamma}\|_C \cdot 2^{-N}
\end{equation*}
where $\|\mathrm ds_{|\Gamma}\|_C$ is the Carleson norm of arclength
on $\Gamma$.

Next we focus our attention on the higher order terms, and give the
estimate for $E_{\mathcal B, 2}$. From~\eqref{eq:boundaryarcs}
and~\eqref{eq:zerocancelling} and the inequality $(a + b)^2 \leq 2(a^2
+ b^2)$ we see that $E_{\mathcal B, 2}$ is bounded by a fixed
multiple of
\begin{multline*}
   \bigl( 1 - |z|^2 \bigr)^2
      \sum_{\Gamma_{i,k} \in \mathcal B} \int_{\Gamma_{i,k}}
      \bigl( 1 - |\xi|^2 \bigr)^2 \biggl| \frac1{|1 - \bar \xi z|^2}
      - \frac1{|1 - \bar \xi_{i,k} z|^2} \biggr|^2
      \, \mathrm d\mu(\xi) \\ + \bigl( 1 - |z|^2 \bigr)^2
      \sum_{\Gamma_{i,k} \in \mathcal B} \int_{\Gamma_{i,k}}
      \frac{(|\xi_{i,k}|^2 - |\xi|^2)^2}{|1 - \bar \xi_{i,k} z|^4} 
      \, \mathrm d\mu(\xi) .
\end{multline*}
For the first term, we use as above the
estimate~\eqref{eq:cancellation} as well as $1 - |\xi| \leq 2^{-2N} (1
- |z|)$ and $|\xi - \xi_{i,k}| \leq 2 \cdot 2^N (1 - |z|)$. Then we
find
\begin{align*}
  \bigl( 1 - &|z|^2 \bigr)^2
      \sum_{\Gamma_{i,k} \in \mathcal B} \int_{\Gamma_{i,k}}
      \bigl( 1 - |\xi|^2 \bigr)^2 \biggl| \frac1{|1 - \bar \xi z|^2}
      - \frac1{|1 - \bar \xi_{i,k} z|^2} \biggr|^2
      \, \mathrm d\mu(\xi) \\
    &\leq 2^4 \cdot 2^{-N} \cdot \bigl( 1 - |z|^2 \bigr) \sum_{\Gamma_{i,k}
      \in \mathcal B} \int_{\Gamma_{i,k}} \bigl( 1 - |\xi|^2 \bigr)
      \biggl| \frac1{|1 - \bar \xi z|^2} - \frac1{|1 - \bar \xi_{i,k}
      z|^2} \biggr| \, \mathrm d\mu(\xi) .
\end{align*}
Observe that the last sum is just~\eqref{eq:movingmodulus} and by the
earlier argument the last expression is bounded by $2^7 \cdot \|\mathrm
ds_{|\Gamma}\|_C \cdot 2^{-2N}$.

For the second term we use that $|1 - \bar \xi_{i,k} z| \geq 1 - |z|$,
$| |\xi_{i,k}| - |\xi| | \leq 2^{-2N} (1 - |z|)$ and $| |\xi_{i,k}| -
|\xi| | \leq \ell(\Gamma_{i,k})$ to arrive at
\begin{align*}
  \bigl( 1 - |z|^2 \bigr)^2
      \sum_{\Gamma_{i,k} \in \mathcal B} \int_{\Gamma_{i,k}}
      &\frac{(|\xi_{i,k}|^2 - |\xi|^2)^2}{|1 - \bar \xi_{i,k} z|^4} 
      \, \mathrm d\mu(\xi) \\
   &\leq 2^4 \cdot \bigl( 1 - |z| \bigr)^{-2}
      \sum_{\Gamma_{i,k} \in \mathcal B} \int_{\Gamma_{i,k}}
      \bigl| |\xi_{i,k}| - |\xi| \bigr|^2 \, \mathrm d\mu(\xi) \\
   &\leq 2^4 \cdot 2^{-2N} \cdot \bigl( 1 - |z|
      \bigr)^{-1} \sum_{\Gamma_{i,k} \in \mathcal B} \ell(\Gamma_{i,k})
   \leq 2^4 \cdot \|\mathrm ds_{|\Gamma}\|_C \cdot 2^{-N} .
\end{align*}
Thus we get that $E_{\mathcal B, 2} \leq C \cdot (2^4 + 1) \cdot
\|\mathrm ds_{|\Gamma}\|_C \cdot 2^{-N}$ for big $N$.

For the short arcs $\Gamma_{i,k} \in \mathcal S_n$, $n \geq N + 1$ we
will use similar estimates as above, however we do not need to be as
delicate. For these arcs, we can use that $|\log x| \leq |1 - x^2|$ to
obtain
\begin{multline*}
  E_{\mathcal S} \overset{\text{def}}{=} \biggl| \sum_{n=N+1}^\infty
      \sum_{\Gamma_{i,k} \in \mathcal S_n}
      \int_{\Gamma_{i,k}} \log \frac{\rho(z, \xi)}{\rho(z, \xi_{i,k})}
      \, \mathrm d\mu(\xi) \biggr| \\
    \leq \sum_{n=N+1}^\infty \sum_{\Gamma_{i,k} \in \mathcal S_n}
      \int_{\Gamma_{i,k}} \bigl| 1 - \frac{\rho(z, \xi)^2}
      {\rho(z, \xi_{i,k})^2} \bigr| \, \mathrm d\mu(\xi) .
\end{multline*}
The same calculations that gave~\eqref{eq:zerocancelling} show that
\begin{equation*}
  \bigl| 1 - \frac{\rho(z, \xi)^2}{\rho(z, \xi_{i,k})^2} \bigr|
    \leq 2 \bigl( 1 - |z|^2 \bigr) \biggl( \biggl|
      \frac{1 - |\xi|^2}{|1 - \bar \xi z|^2} -
      \frac{1 - |\xi|^2}{|1 - \bar \xi_{i,k} z|^2} \biggr| + 
      \frac{\bigl| |\xi_{i,k}|^2 - |\xi|^2 \bigr|}{|1 - \bar \xi_{i,k} z|^2}
      \biggr) .
\end{equation*}
For $\xi \in \Gamma_{i,k} \in \mathcal S_n$, using $|1 - \bar \xi z|
\geq 2^{n-3} (1 - |z|)$ we get
\begin{equation*}
  \bigl( 1 - |\xi|^2 \bigr) \biggl| \frac1{|1 - \bar \xi z|^2} -
      \frac1{|1 - \bar \xi_{i,k} z|^2} \biggr|
    \leq 2^{11} \frac{\bigl( 1 - |\xi| \bigr) \bigl| \xi - \xi_{i,k} \bigr|}
      {2^{3n} \bigl( 1 - |z| \bigr)^3}
    \leq 2^{11} \frac{|\xi - \xi_{i,k}|}{2^{2n} \bigl( 1 - |z| \bigr)^2} .
\end{equation*}
Similarly
\begin{equation*}
  \frac{\bigl| |\xi_{i,k}|^2 - |\xi|^2 \bigr|}{|1 - \bar \xi_{i,k} z|^2}
    \leq 2^7 \frac{|\xi - \xi_{i,k}|}{2^{2n} \bigl( 1 - |z| \bigr)^2} .
\end{equation*}
Adding up, we obtain
\begin{equation*}
  \bigl| 1 - \frac{\rho(z, \xi)^2}{\rho(z, \xi_{i,k})^2} \bigr|
    \leq 2^{14} \frac{|\xi - \xi_{i,k}|}{2^{2n} \bigl( 1 - |z| \bigr)} .
\end{equation*}
Hence
\begin{equation*}
  E_{\mathcal S}
    \leq 2^{14} \sum_{n=N+1}^\infty \frac1{2^{2n} \bigl( 1 - |z| \bigr)}
      \sum_{\Gamma_{i,k} \in \mathcal S_n} \ell(\Gamma_{i,k})
    \leq 2^{14} \cdot \|\mathrm ds_{|\Gamma}\|_C \cdot 2^{-N} .
\end{equation*}

Finally, we estimate the long arcs $\Gamma_{i,k} \in \mathcal L_n$, $n
\geq N+1$. As the zeros on these arcs are well separated, one can
expect only a small contribution from these arcs. We will use an
auxiliary interpolating Blaschke product to find a bound for the
$\mathcal L_n$-terms of~\eqref{eq:decomposable_sum}. By the same
reasoning that led to~\eqref{eq:zerocancelling} and the triangle
inequality,
\begin{align*}
  E_{\mathcal L} &\overset{\text{def}}{=} \biggl| \sum_{n=N+1}^\infty
      \sum_{\Gamma_{i,k} \in \mathcal L_n} \int_{\Gamma_{i,k}} \log
      \frac{\rho(z, \xi)}{\rho(z, \xi_{i,k})} \, \mathrm d\mu(\xi) \biggr| \\
    &\leq 2 \sum_{n=N+1}^\infty \sum_{\Gamma_{i,k} \in \mathcal L_n}
      \int_{\Gamma_{i,k}} \bigl( 1 - |z|^2 \bigr) \biggl( \frac{
      1 - |\xi|^2}{|1 - \bar \xi z|^2} + \frac{1 - |\xi_{i,k}|^2}
      {|1 - \bar \xi_{i,k} z|^2} \biggr) \, \mathrm d\mu(\xi) \\
    &\leq 2^2 \sum_{n=N+1}^\infty \sum_{\Gamma_{i,k} \in \mathcal L_n}
      \max_{\xi \in \overline \Gamma_{i,k}} \frac{\bigl( 1 - |z|^2 \bigr)
      \bigl( 1 - |\xi|^2 \bigr)}{|1 - \bar \xi z|^2} .
\end{align*}
For each $\Gamma_{i,k} \in \mathcal L_n$, let $\zeta_{i,k} \in
\Gamma_{i,k}$ be such that
\begin{equation*}
  \frac{1 - |\zeta_{i,k}|^2}{|1 - \bar \zeta_{i,k} z|^2}
    = \max_{\xi \in \overline \Gamma_{i,k}} \frac{1 - |\xi|^2}
      {|1 - \bar \xi z|^2} ,
\end{equation*}
and define $B_\zeta$ to be the Blaschke product with $\{ \zeta_{i,k}
\}$ as zeros. Now we reorder the summation, and sum with respect to
the placement of the $\zeta_{i,k}$ instead. Then
\begin{equation*}
  E_{\mathcal L} \leq 2^3 \cdot \bigl( 1 - |z| \bigr) \sum_{n=0}^\infty
      \sum_{\zeta_{i,k} \in U_n} \frac{1 - |\zeta_{i,k}|^2}
      {|1 - \bar \zeta_{i,k} z|^2}
\end{equation*}
where $U_0 = Q_z$ and $U_n = 2^n Q_z \setminus 2^{n-1} Q_z$ for $n
\geq 1$. The scaling property~\eqref{eq:scaling} implies that at most
four of the points $\zeta_{i,k}$ are contained in $2^{N-1} Q_z$. These
must be close to the boundary, so that
\begin{equation*}
  2^3 \cdot \bigl( 1 - |z| \bigr) \sum_{n=0}^{N-1} \sum_{\zeta_{i,k} \in U_n}
      \frac{1 - |\zeta_{i,k}|^2}{|1 - \bar \zeta_{i,k} z|^2}
    \leq 4 \cdot 2^4 \cdot 2^{-2N} .
\end{equation*}
For the rest of the terms, we then get
\begin{equation*}
  2^3 \cdot \bigl( 1 - |z| \bigr) \sum_{n=N}^\infty \sum_{\zeta_{i,k} \in U_n}
      \frac{1 - |\zeta_{i,k}|^2}{|1 - \bar \zeta_{i,k} z|^2}
    \leq 2^8 \sum_{n=N}^\infty \frac1{2^n} \sum_{\zeta_{i,k} \in U_n}
      \frac{1 - |\zeta_{i,k}|}{2^n \bigl( 1 - |z| \bigr)}
    \leq 2^9 \cdot C_\zeta \cdot 2^{-N} ,
\end{equation*}
where $C_\zeta$ is the Carleson norm of the measure $\sum (1 -
|\zeta_{i,k}|) \delta_{\zeta_{i,k}}$, which is bounded by a fixed
multiple of $\|\mathrm ds_{|\Gamma}\|_C$. Thus $E_{\mathcal L} \leq
2^9 \cdot (C_\zeta + 1) \cdot 2^{-N}$.

We have now estimated the contribution from all the arcs
$\Gamma_{i,k}$, and we have found that for some constant $C$,
\begin{equation*}
  \bigl| \log |B_1(z)| - \log |I_1(z)| \bigr| \leq C \cdot 2^{-N} .
\end{equation*}
This means that given $\varepsilon > 0$, taking $N$ so that $C \cdot
2^{-N} < \frac\varepsilon2$, we obtain
\begin{equation*}
  \bigl| |B_1(z)| - |I_1(z)| \bigr| < \tfrac\varepsilon2 ,
\end{equation*}
which was what we needed.

\def\cprime{$'$}
\providecommand{\bysame}{\leavevmode\hbox to3em{\hrulefill}\thinspace}


\begin{thebibliography}{10}

\bibitem{Carleson58}
Lennart Carleson,
\emph{An interpolation problem for bounded analytic functions},
Amer.\ J.\ Math. \textbf{80} (1958), 921--930.
\MR{0117349 (22:8129)}

\bibitem{Carleson62}
\bysame,
\emph{Interpolations by bounded analytic functions and the corona problem},
Ann.\ of Math. (2) \textbf{76} (1962), 547--559.
\MR{0141789 (25:5186)}

\bibitem{Frostman35}
Otto Frostman,
\emph{Potentiel d'\'equilibre et capacit\'e des ensembles avec
quelques applications \'a la th\'eorie des fonctions},
Medd.\ Lund.\ Univ.\ Math.\ Sem. \textbf{3} (1935), no.~3, 1--118.

\bibitem{Garnett81}
John~B. Garnett,
\emph{Bounded analytic functions},
Academic Press, 1981.
\MR{0628971 (83g:30037)}

\bibitem{Garnett96}
John~B. Garnett and Artur Nicolau,
\emph{Interpolating {B}laschke products generate {$H^\infty$}},
Pacific J.\ Math. \textbf{173} (1996), no.~2, 501--510.
\MR{1394402 (97f:30050)}

\bibitem{Havin94}
Victor~P. Havin and Nikolai~K. Nikol{\cprime}ski{\u\i} (eds.),
\emph{Linear and complex analysis. {P}roblem book 3. {P}art {II}},
Lecture Notes in Mathematics, vol. 1574, Springer, Berlin, 1994.
\MR{1334345 (96c:00001a)}

\bibitem{Jones81}
Peter~W. Jones,
\emph{Ratios of interpolating {B}laschke products},
Pacific J.\ Math. \textbf{95} (1981), no.~2, 311--321.
\MR{0632189 (82m:30032)}

\bibitem{Lyubarskii01}
Yurii Lyubarskii and Eugenia Malinnikova,
\emph{On approximation of subharmonic functions},
J.\ Anal.\ Math. \textbf{83} (2001), 121--149.
\MR{1828489 (2002b:30043)}

\bibitem{Marshall76}
Donald~E. Marshall,
\emph{Blaschke products generate {$H^\infty$}},
Bull.\ Amer.\ Math.\ Soc. \textbf{82} (1976), 494--496.
\MR{0402054 (53:5877)}

\bibitem{Marshall96}
Donald~E. Marshall and Arne Stray,
\emph{Interpolating {B}laschke products},
Pacific J.\ Math. \textbf{173} (1996), no.~2, 491--499.
\MR{1394401 (97c:30042)}

\bibitem{Mortini04}
Raymond Mortini and Artur Nicolau,
\emph{Frostman shifts of inner functions},
J.\ Anal.\ Math. \textbf{92} (2004), 285--326.
\MR{2072750 (2005e:30088)}

\bibitem{Nicolau05}
Artur Nicolau and Daniel Su\'arez,
\emph{Approximation by invertible functions of {$H^\infty$}},
To appear in Math.\ Scandinavica.

\bibitem{Nikolskii86}
Nikolai~K. Nikol{\cprime}ski{\u\i},
\emph{Treatise on the shift operator},
Springer, Berlin, 1986.
\MR{0827223 (87i:47042)}

\bibitem{Seip04}
Kristian Seip,
\emph{Interpolation and sampling in spaces of analytic functions},
American Mathematical Society, Providence, RI, 2004,
University Lecture Series 33.
\MR{2040080 (2005c:30038)}

\end{thebibliography}
\end{document}